\documentclass[10pt]{article}
\usepackage{amsmath}
\usepackage{amsfonts}
\usepackage{fullpage}
\usepackage{amssymb}
\usepackage{amsthm}
\newcommand{\setlinespacing}[1]%
           {\setlength{\baselineskip}{#1 \defbaselineskip}}

\newcommand{\non}{ ^{*}}
\newcommand{\stan}{^{\circ}\!}

\newcommand{\dt}{\!\vartriangle\!\! t}
\newcommand{\DM}{\textrm{dim}}
\newcommand{\AT}{|A_{\mathbf{T}}|}
\newcommand{\at}{A_{\mathbf{T}}}
\newcommand{\filt}[1]{\langle #1 \rangle_{\mathcal{U}}}
\newcommand{\measb}[1]{\textrm{meas}_{\beta} #1}

\newcommand{\nsn}{ {}^{*}\!\mathbb{N}}
\newcommand{\nsa}{ ^{*}A}
\newcommand{\nsr}{ {}^{*}\mathbb{R}}
\newcommand{\nss}{ ^{*}\mathbb{S}}
\newcommand{\nsx}{ ^{*}\!X}
\newcommand{\st}{\textrm{st}_{\mathbf{T}}^{-1}}
\theoremstyle{plain}
\newtheorem{thm}{Theorem}[section]
\newtheorem{cor}[thm]{Corollary}
\newtheorem{lem}[thm]{Lemma}
\newtheorem{prop}[thm]{Proposition}
\newtheorem{defn}{Definition}[section]
\numberwithin{equation}{subsection}

\begin{document}

\title{NONSTANDARD ANALYSIS, FRACTAL PROPERTIES AND BROWNIAN MOTION}
\author{PAUL POTGIETER}

\date{}

\maketitle \vspace{-1.2cm}
\begin{center}
 \emph{Department of Decision Sciences, University
of South Africa\\ P.O. Box 392, Pretoria 0003, South Africa}\\
\texttt{potgip@unisa.ac.za}
\end{center}

%{\bf{Journal:}} Fractals
 \setcounter{section}{0}

\begin{abstract}
\noindent In this paper I explore a nonstandard formulation of
Hausdorff dimension. By considering an adapted form of the counting
measure formulation of Lebesgue measure, I prove a nonstandard
version of Frostman's lemma and find that Hausdorff dimension can be
computed through a counting argument rather than by taking the
infimum of a sum of certain covers. This formulation is then applied
to obtain a simple proof of the doubling of the dimension of certain
sets under a Brownian motion.

\vspace{0.5cm} \noindent \emph{Keywords: Frostman's lemma,
Nonstandard Hausdorff dimension, Brownian motion}
\end{abstract}

\section{INTRODUCTION}
Using Loeb measure theory, it is possible to construct Lebesgue or
even Wiener measure as a hyperfinite counting measure. In this paper
I explore an extension of the idea to Hausdorff measure, which
yields a hyperfinite formulation of Hausdorff dimension. In cases
where the problem of dimension can be interpreted as an equivalent
problem on a hyperfinite time line, this can lead to a simple and
intuitively satisfying proof. For instance, I shall later consider
certain properties of Brownian motion, and present nonstandard
proofs which are somewhat easier than the original, and also seem to
obey certain  statistical ``rules of thumb".

I now provide a short overview of the necessary nonstandard analysis
as well as the standard formulation of Hausdorff dimension, since
the nonstandard version will follow the same style and notation.

\section{AN INTRODUCTION TO NONSTANDARD ANALYSIS}

Before defining Loeb measures, we briefly introduce the nonstandard
universe in which we will be working. Nonstandard analysis was
introduced in the 1960s by Abraham Robinson~\cite{Robinson}. This
exposition is largely based on the very clear monograph of
Cutland~\cite{Cutland}. Although Loeb measures are standard
measures, their construction involves nonstandard analysis (NSA).
\subsection{The hyperreals} We construct a real line $\nsr$
which is richer than the standard reals $\mathbb{R}$.
This is an ordered field which extends the real numbers to include
non-zero infinitesimals; that is, numbers the absolute value of
which is smaller than any real number; and also positive and
negative ``infinite" numbers. We now make these notions precise.

%We say that $x,y \in\,\!\! \nsr$ are infinitely close whenever $x-y$
%is infinitesimal and denote it by $x\approx y$. Thus, $x\approx y$
%if for every $\varepsilon >0$ in $\mathbb{R}$, $|x-y|<\varepsilon$.
%The set of all such $y$ which are infinitesimally close to $x$ is
%called the \emph{monad} of $x$.\\
There are several ways of constructing the extended universe. Here
we use an ultrapower construction. An axiomatic approach is also
possible, as for instance in~\cite{Robert}. We rather use the
ultrapower because it is pertinent to later constructions, and seems
somewhat more satisfactory in a fractal context.
\begin{defn}
A free ultrafilter $\mathcal{U}$ on $\mathbb{N}$ is a collection of
subsets of $\mathbb{N}$ that is closed under finite intersections
and supersets (i.e. $A\subseteq B$ and $A\in \mathcal{U}$ implies
$B\in \mathcal{U}$), contains no finite sets and for every
$A\subseteq \mathbb{N}$ has either $A\in \mathcal{U}$ or $\mathbb{N}
\setminus A \in \mathcal{U}$.
\end{defn}
Given such a free ultrafilter $\mathcal{U}$ on $\mathbb{N}$ we
construct $\nsr$ as an ultrapower of the reals
\[\nsr = \mathbb{R}^{\mathbb{N}} / \mathcal{U}.\]
The set $\nsr$ that we obtain therefore consists of equivalence
classes of sequences of reals under the equivalence relation
$\equiv_{\mathcal{U}}$, where
\[(a_n )\equiv_{\mathcal{U}} (b_n) \Leftrightarrow \{n: a_n =b_n
\}\in \mathcal{U}.\]  The equivalence class of a sequence $(a_n )$
is denoted by either $(a_n )_{\mathcal{U}}$ or, in the sequel, by
$\langle a_n \rangle_{\mathcal{U}}$. These are what we refer to as
the \emph{hyperreals}. It is clear that $\nsr$ is then an extension
of $\mathbb{R}$, the usual real numbers represented by equivalence
classes of constant sequences.
%We can now also extend subsets of
%$\mathbb{R}$ to $\nsr$; if $A\subset \mathbb{R}$, then $\nsa \subset
%\nsr$ is the

The usual algebraic operations such as $+, \times , <$ are extended
in the above way, but shall be denoted in the usual way. For
example, given $x=\langle x_n \rangle_{\mathcal{U}}$ and $y =
\langle y_n \rangle_{\mathcal{U}}$ in $\nsr$, we can say that $x<y$
if the set of $n\in \mathbb{N}$ on which $x_n < y_n$ is a member of
the ultrafilter ${\mathcal{U}}$.

We usually distinguish three important classes of hyperreals.
Intuitively, the \emph{infinitesimals} are equivalence classes of
sequences converging to $0$ (and therefore smaller than any real
number), \emph{bounded} (or \emph{finite}) hyperreals are
equivalence classes of convergent sequences (including the
infinitesimals), and the \emph{hyperfinite} or \emph{infinite}
hyperreals are equivalence classes of sequences diverging to
$\infty$ or $-\infty$ (and thus, in absolute value, larger than any
real number).

There exists a function (which we'll discuss in more detail when
considering nonstandard topology) from the finite hyperreals to
$\mathbb{R}$ which associates to each nonstandard number an element
of $\mathbb{R}$, called the \emph{standard part} of the number, as
expressed in the following theorem:
\begin{thm}
If $x\in \!  \nsr$ is finite (that is, $-R< x < R$ for some standard
$R\in \mathbb{R}$), then there is a unique $r\in \mathbb{R}$ such
that $x\approx r$. Any finite hyperreal is thus expressible as
$x=r+\delta$ with $r\in \mathbb{R}$ and $\delta$ infinitesimal.
\end{thm}
We call $r$ in the above theorem the \emph{standard part} of $x$ and
denote it as either $ ^{\circ} x$ or as $\textrm{st} (x)$. Both are
used, sometimes in conjunction, to improve readability.

For infinitesimals the standard part is clearly $0$; indeed, this
suffices as a definition of infinitesimal for nonzero numbers. We
say that $x,y \in \! \nsr$ are infinitely close whenever $x-y$ is
infinitesimal and denote it by $x\approx y$. Thus, $x\approx y$ if
for every $\varepsilon >0$ in $\mathbb{R}$, $|x-y|<\varepsilon$. The
set of all such $y$ which are infinitesimally close to $x$ is called
the \emph{monad} of $x$.

Functions are also defined by pointwise operations. Given a function
$f: \mathbb{R}\to \mathbb{R}$, a nonstandard counterpart of $f$ is
given by the function $F: \nsr \to \nsr$, defined by
\[F(\langle a_n \rangle_{\mathcal{U}}) = \langle f(a_n
)\rangle_{\mathcal{U}}.\] Note that this is applicable to the
characteristic functions of sets as well, giving us a way to extend
$A\subset \mathbb{R}$ to its nonstandard counterpart $\nsa \in
\nsr$. (The imbedding $ ^{*}:\mathbb{R}\to \nsr$ so obtained is a
Boolean homomorphism.) If $f$ is a real function defined on a set
$A$, we call the function $F$ defined on $^{*}\! A$ by the above a
\emph{lifting} of $f$, and denote it by $^{*}\! f$. Note that if
$t\in \mathbb{R}$, then $f(t) = \stan \! F(t)$, and if $\tau \in \!
\nsr$ then $^{*}\! f(\tau) \approx F(\tau)$.

 Exactly which properties
of $\nsr$ are inherited from $\mathbb{R}$ is specified in the
following theorem, a restricted version of the more general
\emph{transfer principle}:
\begin{thm}
Let $\varphi$ be any first order statement. Then $\varphi$ holds in
$\mathbb{R}$ if and only if $ ^{*} \varphi$ holds in $\!\nsr$.
\end{thm}
A first order statement $\varphi$ (or $ ^{*}\varphi$ in $\nsr$) is
one referring to elements (fixed or variable) of $\mathbb{R}$
(respectively, $\!\nsr$) and to fixed functions and relations on
$\mathbb{R}$ (respectively, $\nsr$), that uses the usual logical
connectives \emph{and} ($\wedge$), \emph{or} ($\vee$),
\emph{implies} ($\rightarrow$) and \emph{not} ($\neg$).
Quantification may be done over elements but not over relations or
functions; i.e., $\forall x$, $\exists y$ are allowed, but $\forall
f$, $\exists R$ are not. As an example, the density of the rationals
in the reals can be written as
\[ \forall x \forall y (x<y \rightarrow \exists z (z\in \mathbb{Q}
\wedge (x<z<y))),\] an expression meaning, ``between every two reals
is a rational". From the transfer principle we can therefore
immediately conclude that the statement is true in $\nsr$, i.e. that
the hyperrationals are dense in the hyperreals.  The following
theorem relates the concepts of convergence and being infinitely
close.
\begin{thm}
Let $(s_n )$ be a sequence of real numbers and let $l\in
\mathbb{R}$. Then
\[s_n \to l \textrm{ as } n\to \infty \Longleftrightarrow\,\!
^{*}\! s_K \approx l \textrm{ for all infinite } K\in\,\! \nsn.\]
\end{thm}
{\bf{Proof.}}~\cite{Cutland} Suppose that $s_n \to l$ and let $K\in
\nsn$ be a fixed infinite number. We must show, for all real
$\varepsilon >0$, that $| ^{*}\! s_K -l|<\varepsilon$. From ordinary
real analysis we know that there exists some $n_0 \in \mathbb{N}$
such that
\[\forall n\in \mathbb{N} \, [n\geq n_0 \rightarrow |s_n
-l|<\varepsilon].\] According to the transfer principle, the
following is true in $\nsr$:
\[\forall N\in\,\! \nsn \, [N\geq n_0 \rightarrow | ^{*}\!s_N -l|\leq
\varepsilon].\]
In particular, $| ^{*}\!s_K -l|<\varepsilon$ as required.\\

Conversely, suppose that $ ^{*}\!s_K \approx l$ for all infinite
$K\in\,\! \nsn$. For any given real $\varepsilon >0$ we have
\[\exists K\in\,\! \nsn \, \forall N\in\,\! \nsn \, [N\geq K \rightarrow |
^{*}\!s_N -l|<\varepsilon].\] By transferring this ``down" to
$\mathbb{R}$, we get
\[\exists k\in \mathbb{N}\,\ \forall n\in \mathbb{N} \, [n\geq k
\rightarrow |s_n -l|<\varepsilon].\] This implies convergence to
$l$.\hfill{$\square$}

\subsection{The nonstandard universe}
The principles of the previous section can be used in a much broader
context than just real analysis. Given any mathematical object
$\mathcal{M}$ (whether it is a group, ring, vector space, etc.), we
can construct a nonstandard version $ ^{*}\!\mathcal{M}$. We
consider a somewhat more economical construction however, by
starting with a working portion of the mathematical universe
$\mathbb{S}$ and ending up with a $ ^{*}\mathbb{S}$ which will
contain $ ^{*}\!\mathcal{M}$ for every $\mathcal{M} \in \mathbb{S}$.
This has the added advantage of preserving some of the relations
between structures through the more general transfer principle.\\

We start with the superstructure over $\mathbb{R}$, denoted by
$\mathbb{S}=S (\mathbb{R})$. It is defined as follows:
\begin{eqnarray*}
S_{0}(\mathbb{R}) &=& \mathbb{R}\\
S_{n+1}(\mathbb{R}) &=& S_n (R) \cup \mathcal{P}(S_n
(\mathbb{R})), \quad n\in \mathbb{N}\\
\mathbb{S} &=&\bigcup_{n\in \mathbb{N}} S_n ({\mathbb{R}}).
\end{eqnarray*}
($\mathcal{P}(A)$ denotes the power set of the set $A$.)\\

If a larger (or simply different) universe is required, start the
same process with a more suitable set than $\mathbb{R}$. However,
one usually only needs the first few levels of this construction.

Next one must construct a mapping $ ^{*} :S(\mathbb{R}) \to S(\nsr)$
associating to an $\mathcal{M}\in \mathbb{S}$ a nonstandard
extension $ ^{*}\!\mathcal{M} \in S(\nsr)$. The nonstandard universe
can now be constructed by means of an ultrapower
\[\mathbb{S}^{\mathbb{N}} /\mathcal{U}.\]
This is somewhat more complicated to do than in the case of $\nsr$
and we do not go into details here. It is sufficient to consider the
nonstandard universe as the set of objects
\[ ^{*}\mathbb{S} = \{x:x\in\,\!\! ^{*}\mathcal{M} \textrm{ for some }
\mathcal{M} \in \mathbb{S}\}.\]

Sets in $\nss$ are called \emph{internal} sets. It should be noted
that $\nss \in S(\nsr)$, but that $S(\nsr)$ contains sets which
are not internal.\\

We now also have a Transfer Principle which specifies which
statements may be moved from one structure to the other (see
\cite{Albeverio}, for instance). A \emph{bounded quantifier
statement} differs from a first-order statement only in that
quantifiers must range over a fixed set - which may mean a set of
functions or relations, not just elements of $\mathbb{S}$ (the
functions and relations are of course merely elements of the
superstructure). Thus, quantifiers like $\forall x\in A$ or $\exists
y\in B$ are allowed, but not unbounded quantifiers such as $\forall
x$ and $\exists y$. Note that boundedness of the quantifier is often
implied in the exposition and is not always specifically indicated
in the statement.
\begin{thm}
A bounded quantifier statement $\varphi$ holds in $\mathbb{S}$ if
and only if $ ^{*}\varphi$ holds in $\nss$.
\end{thm}

We show as a first application that the concept of supremum
transfers:

\begin{prop}
Every nonempty internal subset of $\nsr$ with an upper bound has a
least upper bound.
\end{prop}
\textbf{Proof.}~\cite{Cutland} The notation used in this proof
refers back to our construction of the nonstandard universe. We
express the fact that any nonempty internal subset has a least upper
bound by the statement
\begin{eqnarray*}\Phi (\nsr ,\non\! S_2 (\mathbb{R})) =\forall A\in
\non\! S_2 (\mathbb{R})[A\neq \emptyset \wedge
(\exists x\in\,\!\! \nsr\, (\forall y\in A(y<x)))\rightarrow&& \\
\exists z\in\,\!\! \nsr \, (\forall y \in A (y<x) \wedge \forall
u\in\,\!\! \nsr \, \forall y\in A\,(y\leq u\rightarrow z\leq u))].&&
\end{eqnarray*} Since the condition $\Phi (\mathbb{R}, S_2
(\mathbb{R}))$ is true in $S(\mathbb{R} )$, the transferred
condition is true in $S(\nsr )$.\hfill{$\square$}

The Transfer Principle can after some consideration be seen to apply
only to \emph{internal} sets. For instance, the concept of supremum
implies that each bounded set will have a least upper bound.
However, $\mathbb{N}$ seen as a member of $\nsr$ is bounded by each
element of $\nsn \setminus \mathbb{N}$, but has no supremum. It is
therefore an \emph{external} (i.e., non-internal) set. It is
sufficient to think of internal sets as sets in $^{*}\!\mathbb{S}$
that are images under the homomorphism $^{*}$.

Another important property of any nonstandard universe constructed
as an ultrapower is that of $\aleph_1$-saturation:
\begin{defn}
If $(A_m)_{m\in \mathbb{N}}$ is a countable decreasing sequence of
nonempty internal sets, then $\cap_{m\in \mathbb{N}} A_m \neq
\emptyset.$
\end{defn}
A useful reformulation of this is known as \emph{countable
comprehension}: Given any sequence $(A_n )_{n\in \mathbb{N}}$ of
internal subsets of an internal set $A$, there is an \emph{internal}
sequence $(A_n )_{n\in \nsn}$ of subsets of $A$ that extends the
original sequence (that is, agrees with the original sequence on the
standard natural numbers, and inherits any properties that can be
phrased as a bounded quantifier statement). This property is used in
the construction of Loeb measure. Two other important properties are
those of \emph{overflow} and \emph{underflow}: the first states that
an internal set containing arbitrarily large finite numbers must
contain an infinite number, and the second that an internal set
containing arbitrarily small infinite numbers must contain a finite
number. (Considering the reciprocals of the numbers, one can state
similar results for infinitesimals.)

We observe as well that the notion of the cardinality of a set
transfers. For a finite set the cardinality function simply gives
the number of elements of a set (more formally and generally, the
cardinality of a set $X$ is the least ordinal $\alpha$ such that
there is a bijection between $X$ and $\alpha$). In the sets I will
consider in the case of nonstandard Hausdorff dimension, all will be
finite within the nonstandard context. One may therefore still
intuitively regard the transferred cardinality function as giving
the ``size" of a set. The cardinality of a set $X$ will be denoted
by $|X|$ in both standard and nonstandard cases; which is meant
should be clear.

\subsection{Nonstandard topology}
Since Brownian motion, and hence continuous functions, is later
considered, a knowledge of nonstandard topology is required.
 Firstly, we see that the concept of being
infinitely close, and therefore the idea of a \emph{monad}, can be
extended:
\begin{defn}
Let $(X,\tau)$ be a topological space.
\begin{enumerate}
\item[(i)]{For $a\in X$ the monad of $a$ is
\[monad(a)= \bigcap_{a\in U\in \tau }\!  {}^{*}\!U.\]}
\item[(ii)]{For $x\in \,\nsx$, we write $x\approx a$ if $x\in
\textrm{monad}(a)$.}
\item[(iii)]{ $x\in \,\nsx$ is said to be nearstandard if $x\approx a$ for some $a \in X$.}
\item[(iv)]{For any $Y\subseteq \,\nsx$, we denote the
nearstandard points in $Y$ by $\textrm{ns}(Y)$.}
\item[(v)]{$\textrm{st}(Y)= \{a\in X: x\approx a \textrm{ for some } x\in
Y\}$ is called the standard part of $Y$ (also denoted by $\stan
Y$).}
\end{enumerate}
\end{defn}
The following result allows us to generalise the pointwise standard
part mapping:
\begin{prop}
A topological space $X$ is Hausdorff if and only if
\[monad(a) \cap \textrm{monad}(b) = \emptyset \textrm{
for } a\neq b, \quad a,b\in X.\]
\end{prop}
This means we can define the function
\[\textrm{st}: \textrm{ns}(\non X)\to X\]
as
\[\textrm{st}(x) = \textrm{ the unique }a\in X \textrm{ with }
a\approx x.\] The notation $\stan x= \textrm{st}(x)$
is again used interchangeably.\\

We mention some general topological results.
\begin{prop}~\cite{Cutland}
Let $(X,\tau )$ be separable and Hausdorff. Suppose $Y\subseteq
\non\!X$ is internal and $A\subseteq X$. Then
\begin{enumerate}
\item[(i)]{$\textrm{st}(Y)$ is closed,}
\item[(ii)]{if $X$ is regular and $Y\subseteq \textrm{ns}(\non
X)$, then $\textrm{st}(Y)$ is compact,}
\item[(iii)]{$\textrm{st}(\non A) = \overline{A}$ (closure of
$A$),}
\item[(iv)]{if $X$ is regular, then $A$ is relatively compact iff
$\non A \subseteq \textrm{ns} (\non X)$.}
\end{enumerate}
\end{prop}
Since we will be dealing almost exclusively with continuous
functions, we should introduce corresponding notions in the
nonstandard universe.
\begin{defn}
Let $Y$ be a subset of $\nsx$ for some topological space $X$  and
let $F: \,\nsx \to \nsr$ be an internal function. Then $F$ is said
to be \emph{S-continuous on Y} if for all $x,y \in Y$ we have
\[x \approx y \Rightarrow F(x)\approx F(y).\]
\end{defn}
The following result allows us to switch from the one notion of
continuity to another.
\begin{thm}~\cite{Cutland}
If $F:\nsr \to \nsr$ is $S$-continuous on an interval $^{*}\![a,b]$
for real $a,b$ and $F(x)$ is finite for some $x\in ^{*}\![a,b]$,
then the standard function defined in $[a,b]$ by
\[f(t) = \stan F(t)\]
is continuous and $ ^{*}\! f(\tau) \approx F(\tau)$ for all $\tau
\in ^{*}\![a,b]$.
\end{thm}
%Given a real function $f$ defined on an interval $[a,b]$, we shall
%call any function $F$ on $^{*}\![a,b]$ such that $f(t)=\stan F(t)$,
%a \emph{lifting} of $f$.

\subsection{Loeb measure}
A Loeb measure is a standard measure, but constructed from a
nonstandard one. That is, the Loeb measure exists on a
$\sigma$-algebra and obeys all the usual rules for a measure, e.g.
countable additivity.

We start with a given internal set $\Omega$ and an algebra
$\mathcal{A}$ of internal subsets of $\Omega$. Let $\mu$ be a
finitely additive finite internal measure on $\mathcal{A}$. Thus
$\mu$ is a function from $\mathcal{A}$ to $ \non [0,\infty)$ such
that $\mu(\Omega)<\infty$ and $\mu(A\cup B)=\mu(A)+\mu(B)$ for
disjoint $A,B\in \mathcal{A}$. (We focus only on bounded Loeb
measures; infinite ones shall not concern us in the sequel). We can
then define the mapping
\[\stan \mu: \mathcal{A} \to [0,\infty)\]
by $\stan \mu(A) = \stan\!(\mu(A))$. This is finitely additive and
therefore $(\Omega ,\mathcal{A} ,\stan\!\mu )$ is a standard
finitely additive measure space. This is usually not a measure,
since $\stan \mu$ is usually not $\sigma$-additive. We shall see
shortly, however, that it is \emph{almost} a measure. The following
crucial theorem was proved by Loeb~\cite{Loeb}.
\begin{thm}
There is a unique $\sigma$-additive extension of $\stan \mu$ to the
$\sigma$-algebra $\sigma (A)$ generated by $\mathcal{A}$. The
completion of this measure is the Loeb measure associated with
$\mu$, denoted by $\mu_L$. The completion of $\sigma(A)$ is the Loeb
$\sigma$-algebra, denoted by $L(\mathcal{A})$.
\end{thm}
The more straightforward proof depends on the notion of a Loeb null
set:
\begin{defn}
Let $B\subseteq \Omega$, where $B$ is not necessarily internal. We
call $B$ a \emph{Loeb null set} if for each real $\varepsilon >0$
there is a set $A\in \mathcal{A}$ with $B\subseteq A$ and
$\mu(A)<\varepsilon$.
\end{defn}
This allows us to make precise the notion that $\mathcal{A}$ is
almost a $\sigma$-algebra.:
\begin{lem}~\cite{Cutland}
Let $(A_n )_{n\in \mathbb{N}}$ be an increasing family of sets, with
each $A_n \in \mathcal{A}$ and let $B=\bigcup_{n\in \mathbb{N}}
A_n$. Then there is a set $A\in \mathcal{A}$ such that
\begin{enumerate}
\item[(i)]${B\subseteq A}$ \item[(ii)]{$\stan \mu (A) = \lim_{n\to
\infty} {\,\stan\mu (A_n )}$ and} \item[(iii)]{$A\setminus B$ is
null}
\end{enumerate}
\end{lem}
\textbf{Proof.} Let $\alpha = \lim_{n\to \infty} \,\stan \mu(A_n )$.
For any finite $n$,
\[\mu (A_n )\leq \stan \mu (A_n ) +\frac{1}{n} \leq \alpha
+\frac{1}{n}.\] Let $(A_n )_{n\in \nsn} $ be a sequence of sets in
$\mathcal{A}$ extending the sequence $(A_n )_{n\in \mathbb{N}}$,
possible by $\aleph_1$ saturation. The overflow principle then
guarantees an infinite $N$ such that
\[\mu (A_n )\leq \alpha + \frac{1}{N}.\]
If we now let $A=A_N $, (i) will hold because $A_n \subseteq A$ for
each $n$. Also, $\mu (A_n )\leq \mu (A)$ for finite $n$, so $\stan
\mu (A_n )\leq \mu (A)\leq \alpha$ and therefore $\stan \mu (A)
=\alpha $. This gives (ii). For (iii), note that $A\setminus B
\subseteq A\setminus A_n$ and $\stan \mu (A\setminus A_n )=\stan \mu
(A) -\stan (A_n )\to 0.$ \hfill{$\Box$}

Thus $\mathcal{A}$ is a $\sigma$-algebra modulo null sets. We can
now define the concepts \emph{Loeb measurable} and \emph{Loeb
measure} exactly:
\begin{defn}
%\begin{enumerate}

(i) Let $B\subseteq \Omega $. We say that $B$ is \emph{Loeb
measurable} if there is a set $A\in \mathcal{A}$ such that
$A\bigtriangleup B$ (the symmetric difference of $A$ and $B$) is
Loeb null. The collection of all the Loeb measurable sets is denoted
by $L(\mathcal{A})$. $L(\mathcal{A})$ is known as the \emph{Loeb
algebra}.

(ii) For $B\in L(\mathcal{A})$ define
\[\mu_L (B)=\stan\! \mu (A)\]
for any $A\in \mathcal{A}$ with $A\bigtriangleup B$ null. We call
$\mu_L (B)$ the \emph{Loeb measure} of $B$.
%\end{enumerate}
\end{defn}
This brings us to the central theorem of Loeb measure
theory~\cite{Loeb}.
\begin{thm}
$L(\mathcal{A})$ is a $\sigma$-algebra and $\mu_L$ is a complete
$\sigma$-additive measure on $L(A)$.
\end{thm}
The measure space $\mathbf{\Omega} = (\Omega , L(\mathcal{A}),
\mu_L)$ is called the Loeb space associated with $(\Omega ,
\mathcal{A}, \mu)$. If $\mu (\Omega )=1$, we refer to
$\mathbf{\Omega}$ as a \emph{Loeb probability space}.

\subsection{Loeb counting measure}
I devote a short but separate section to the idea of counting
measures, since they are prominent throughout the sequel, being
useful in the construction of nonstandard Hausdorff dimension.

Let $\Omega = \{1,2,\dots ,N\}$, where $N\in \nsn \setminus
\mathbb{N}$. The set $\Omega$ is internal. Define the counting
probability $\nu$ on $\Omega$ by
\[\nu (A) = \frac{|A|}{N},\]
for $A\in \non \mathcal{P}(\Omega )=\mathcal{A}$. The cardinality
function $|\cdot |$ transfers, so $|A|$ can be interpreted as
denoting the number of elements in $A$. The Loeb counting measure
$\nu_L$ is the completion of the extension to $\sigma (A)$ of the
finitely additive measure $\stan \nu$.

The following definition will be used repeatedly in the sequel:
\begin{defn}
Fix $N\in \nsn \setminus \mathbb{N}$ and let $\dt =N^{-1}$. The
hyperfinite time line for the interval $[0,1]$ based on the
infinitesimal $\dt$ is the set
\[\mathbf{T}=\{0,\dt ,2\!\dt ,3\!\dt ,\dots ,1-\!\dt \}.\]
\end{defn}
The following theorem provides an intuitive construction of Lebesgue
measure.
\begin{thm}
Let $\nu_L$ be the Loeb counting measure on $\mathbf{T}$. Define
\begin{enumerate}
\item[(i)]{$\mathcal{M}=\{B\subseteq [0,1]: \st (B) \textrm{ is
Loeb measurable}\}$, where $\st (B)=\{\textsf{t}\in \mathbf{T}
:\stan \textsf{t}\in B\}$}
\item[(ii)]{$\lambda (B) = \nu_L (\st
(B) )$ for $B\in \mathcal{M}$}
\end{enumerate}
Then $\mathcal{M}$ is the completion of the Borel sets
$\mathcal{B}[0,1]$ and $\lambda$ is Lebesgue measure on
$\mathcal{M}$.
\end{thm}

A sketch of the proof can be found in~\cite{Cutland}

\section{HAUSDORFF DIMENSION}

Given a compact set $A$ on the unit interval (or any bounded subset
of $\mathbb{R}$) and $\epsilon>0$, consider all coverings of the set
by open balls $B_n$ of diameter smaller than or equal to $\epsilon$.
For each cover, form the sum
\[ \sum_{n=0}^{\infty} \|B_n \|^{\alpha} ,\]
where $\|\cdot \|$ denotes the diameter of a set (i.e., the maximum
distance between any two points of the set). For each $A$ and
$\epsilon
>0$, take the infimum over all such sums, as $\{B_n \}$ ranges over
all possible covers of $A$ of diameter $\leq \epsilon$:
\[S_{\alpha}^{\epsilon}(A) = \inf_{\{B_n \}} \sum_{n} \|B_n
\|^{\alpha}.\] As $\epsilon$ decreases to $0$,
$S_{\alpha}^{\epsilon}(B)$ increases to a limit
$\textrm{meas}_{\alpha} (A)$ (which might be infinite) which is
called the $\alpha$-Hausdorff measure of $A$, or the Hausdorff
measure of $A$ in dimension $\alpha$ (I will refer to this as just
``the measure" when the context is clear). Since
$\textrm{meas}_{\alpha}$ is $\sigma$-subadditive but otherwise
satisfies the requirements of a measure, it is an outer measure.

\begin{defn} The \emph{Hausdorff dimension}, $dim A$, of a compact set
$A\subseteq [0,1]$ is the supremum of all the $\alpha \in [0,1]$ for
which, for any cover $B$ of $A$, $\textrm{meas}_{\alpha}(B)=\infty$.
This is equal to the infimum of all $\beta \in [0,1]$ for which
there exists a cover $C$ of $A$ such that
$\textrm{meas}_{\alpha}(B)=0$.
\end{defn}
To see that the supremum of the one set of values is indeed equal to
the infimum of the other, let $0< \alpha <\beta \leq 1$ and consider
the following: \[ \sum_n \|B_n\|^\beta \leq \sup_n \|B_n\|^{\beta
-\alpha} \sum_n \|B_n\|^{\alpha }.\] Hence, if $\textrm{meas}_\alpha
(A) < \infty$, $\textrm{meas}_\beta (A)=0$, or equivalently,
$\textrm{meas}_\alpha (A) =\infty$ if $\textrm{meas}_\beta (A)>0$.

From Hausdorff's original paper~\cite{Hausdorff1} it may inferred
that his intention was somewhat akin to some of the motivation
behind the creation of nonstandard analysis (which shall soon be
using in this context). In this paper, he states:
\begin{quote} In this way, the dimension becomes a sort of characteristic measure of
graduality similar to the `order' of convergence to zero, the
`strength' of convergence, and related concepts.\end{quote}

Although I work  almost exclusively with compact sets in one
(topological) dimension, it is possible to do so in any number of
dimensions. The principles remain inviolate and the Hausdorff
dimension of a set is the same whether we consider it as a subset of
$\mathbb{R}$ or $\mathbb{R}^n$.

\section{NONSTANDARD HAUSDORFF DIMENSION}
In this section and those following I show that a formulation of
Hausdorff measure as a nonstandard counting measure, similar to
Loeb's formulation of Lebesgue measure, is possible and prove some
well-known theorems using these nonstandard techniques. It turns out
that some interesting dimensional properties of Brownian paths
become quite easy to prove using hyperfinite counting arguments.

Before we start the proof, we need a nonstandard version of the
following result~\cite{Frostman}. For notational convenience, the
diameter of a set is denoted by $\|\cdot \|$ and the finite
cardinality function or its transfer by $| \cdot|$, although which
is intended should be clear. The following definition is included
for clarity.

\begin{defn}
A Radon measure is a measure defined on a $\sigma$-algebra of Borel
sets of a set $X$ which is both inner regular (the measure of a set
is the supremum of the measures of the compact sets contained
therein) and locally finite (every point has a neighbourhood of
finite measure).
\end{defn}

\begin{thm}
(Frostman's lemma) Let $A$ be a compact subset of $[0,1]$ and $\beta
\in (0,1)$. Then $\measb{A} >0$ if and only if there exists a
nonzero Radon measure $\mu$ on $A$ such that $\mu (B) \leq C
\|B\|^{\beta}$ for each interval $B \subseteq [0,1]$ and some
positive $C$.
\end{thm}

I will now prove a nonstandard analogue of Frostman's lemma, in
which I use the hyperfinite time line, as introduced in Section 2.5.
Instead of considering the time line as a set of points in $\non
[0,1]$, it will be useful to regard it instead as a subdivision of
the nonstandard real line between $0$ and $1$ into equal parts of
infinitesimal size. The ``time line based on the number $N$"
indicates that these intervals have length $N^{-1}$. If we base the
time line on a sequence $\{a_n \}$, this implies that $\non [0,1]$
is divided into $N=\langle a_n \rangle_{\mathcal{U}}$ equal
intervals of the form $[k\dt, (k+1)\dt]$, $0\leq k <N$.

(Note that I abuse the notation slightly in the following by using $
\stan\left(\frac{|A'|}{2^{N\alpha}}\right)
>0$ to mean either that the standard part of the expression in
brackets exists and is larger than $0$, or that the expression is
infinite.)
\begin{thm}
Let $A$ be a compact subset of $[0,1]$. Suppose $\mathbf{T}$ is a
hyperfinite time line on $[0,1]$, based on the dyadic sequence
$\{2^n \}_{n\geq 1}$, and $A'$ is any internal subset of
$\mathbf{T}$, such that its standard part is $A$. If
 $ \stan\left(\frac{|A'|}{2^{N\alpha}}\right)
>0,$  there exists a nonstandard
measure $\mu$ on $\mathbf{T}$ such that the (nonzero) Loeb measure
$\mu_L$ associated to $\mu$ has the property that for an absolute
constant $C$ and an arbitrary interval $B\subseteq [0,1]$, it is
true that $\mu_L (B)\leq C\| B \|^{\alpha}$.
\end{thm}
\textbf{Proof.} The measure in question is not quite as simple as,
for instance, the counting measure used to generate Lebesgue
measure. In this case we have to take into account how ``close"
elements of $A$ are to each other and a uniform counting measure
cannot provide that information. Thus the construction of the
measure is not generic but will depend specifically on the nature of
$A$.

We use a time line $\mathbf{T}$ based on the hyperfinite number
$2^N$, where $N = \langle 1,2,3,... \rangle_\mathcal{U}$. We say a
dyadic interval is of order $m$ if it has length $2^{-m}$. Let $A'$
be any internal subset of $\mathbf{T}$ such that $\stan A' = A$. On
all intervals of order $N$, we distribute the mass $2^{-N\alpha}$ if
the interval is also in $A'$. Clearly the desired inequality holds
trivially for intervals of order $N$, but may not hold for any
intervals of lower order. Because of the condition on the
cardinality of $A'$, we also have a total mass which is larger than
some positive (standard) real number. We use this mass to normalise,
dividing each interval of $A'$ by this total. Hence the inequality
continues to hold. We now consider order $N-1$. If the inequality
continues to hold on such an interval, we leave it be, and it
retains its original measure. If now, however, two intervals of
order $N$ are both contained in an interval of order $N-1$, the
inequality will be violated, and we need to multiply each by a
factor of $2^\alpha /2$. Once this is accomplished for all intervals
of order $N-1$, we have that the inequality is satisfied for all
these intervals and for those of order $N$ (since their individual
masses cannot increase during this procedure) and again obtain a
total measure larger than some standard positive number. We
normalise the measure thus obtained, and move on to the next level.
At each subsequent level we refine the masses all the way up to the
$N$th level, so that the inequality will continue to hold. Since we
adjust by at most a finite factor at each stage, and the total mass
is larger than some real number, we never have to normalise by
dividing by an infinitesimal. After $N$ steps, we obtain a finite
internal finitely additive nonstandard measure $\mu$ supported by
$A'$ for which there is an absolute (real) constant $C$ such that
the inequality $\mu (B)\leq C\|B\|^{\alpha}$ holds for dyadic
intervals of any order.

(If we only have that $|A'|/2^{N\alpha}\approx 0$, this construction
does not hold. Indeed, in the proofs of Theorems 4.2 and 4.3 it is
shown (independently from this result) that in such a case no such
probability will exist. This construction will fail in such a case
because we require the number $C$ in the statement of the theorem to
be a standard real. If we normalise by dividing by an infinitesimal,
the inequality will not hold for any real number at order $N$. If we
normalise by dividing by an infinite number there is no such
objection, since it will then hold for any real number $C$.)

In the same way as with Lebesgue measure, we now obtain a measure on
the Borel subsets of $[0,1]$ by taking the standard part of $\mu$
and performing the necessary completions, obtaining the measure
$\mu_L$. The inequality on dyadic intervals in $[0,1]$ will then
also hold for $\mu_L$. An arbitrary interval $D$ will always be
contained in two such dyadic intervals and therefore
\[\mu_L (D) \leq C\| D\|^{\alpha}.\]
\hfill{$\square$}

We prove the main result of this section in two separate theorems.
The first guarantees the existence of a subset of a time line from
which we can compute the dimension and the second shows that the
choice of set is not very important. It is proved for subsets of
$[0,1]$ only, but note that it can easily be extended to any compact
interval and arbitrary (finite) dimension.
\begin{thm}
Given a compact subset $A$ of $[0,1]$, there is a subset
$A_{\mathbf{T}}$ on the hyperfinite time line $\mathbf{T}$ and a
hyperfinite number $N\in {}\nsn \setminus \mathbb{N}$ such that
$\stan A_{\mathbf{T}} = A$ and
\begin{eqnarray*}
\stan \left( \frac{|A_{\mathbf{T}}|}{N^{\beta}} \right) &=&
\infty \textrm{ for } \beta<\alpha\\
\stan\left( \frac{|A_{\mathbf{T}}|}{N^{\beta}} \right) &=& 0
\textrm{ for } \beta>\alpha
\end{eqnarray*}
if and only if $dim A =\alpha$.
\end{thm}
{\bf{Proof.}} Suppose that $\beta< \textrm{dim} A$. We know that the
sum diverges to infinity as the sizes of the intervals decrease.
Thus there will be some $N\in \mathbb{N}$ such that the
$\beta$-Hausdorff sum will be larger than $1$ for covers
constituting of sets with diameter smaller than $2^{-K}$, for all
$K>N$.

We will now state, as a bounded quantifier statement, that this will
hold for any cover and that such a cover always exists, a seemingly
trivial point in the standard case, but not as obvious in the
nonstandard.

Let $S=S(A,X,K,J)$ be the following statement, where $X\subset
\mathbb{N}\times \mathbb{N}$:
\[ \begin{array}{lc}
 S = & \forall x\in
A\exists (i,j)\in X \left( x\in
\left(\frac{i}{2^K},\frac{j}{2^K}\right]\right) \wedge \forall
(i,j)\in X \exists x\in A \left( x\in
\left(\frac{i}{2^K},\frac{j}{2^K} \right]\right) \\& \\ & \wedge
\forall (i,j)\in X \left( 2^{-K}\leq (j-i)2^{-K} \leq
2^{-J}\right)\\ & \\ & \wedge \left[ (i,j)\in X \Rightarrow
\neg\left( \exists k\in \{0,1,\dots ,2^K -1\}((j,k)\in X))\right)
\right].\end{array}\] The statement $S$ encapsulates the idea that
there is a cover of $A$ by intervals no smaller than $2^{-K}$, such
that no member of the cover is redundant (i.e., does not contain a
member of $A$). What is more, $S$ states that the largest interval
is no larger than $2^{-J}$, which we may assume because $A$ is
compact, and that no two intervals border in each other, because
then the inequality $2\cdot 2^{-K\beta} > (2\cdot 2^{-K})^{\beta}$
implies that the Hausdorff sum may be decreased.

Let $T=T(X,K,\beta )$ be the statement
\[ \sum_{(i,j)\in X} \left(\frac{j-i}{2^K}\right)^{\beta} >1,\]
which captures the idea that since the Hausdorff sum diverges, it
will eventually be larger than $1$.

We then express $\beta< \textrm{dim} A$ as:
\[\begin{array}{cc} & [\exists N \in \mathbb{N} \forall J> N \forall K\geq J \forall X\subseteq
\{0,1,\dots ,2^K -1\} \times \{0,1,\dots ,2^K -1\}\\ &
(S\Rightarrow T)]\quad \wedge \\
&[\exists N \in \mathbb{N} \forall J> N  \forall K\geq J \exists
X\subseteq \{0,1,\dots ,2^K -1\}\times \{0,1,\dots ,2^K -1\}\\&
(S\Rightarrow T)].
\end{array}\]

The first part of the above statement states that any cover (indexed
by the set $X\subseteq \mathbb{N}\times \mathbb{N}$) of $A$ which
satisfies the conditions stipulated in $S$ will yield a Hausdorff
sum larger than $1$, and the second part states that such a cover
will indeed exist.

The transferred statement now reads as
\[\begin{array}{cc} &[\exists N \in {}^{*}\! \mathbb{N}\forall J>N\forall K\geq J
\forall X\subseteq \{0,1,\dots ,2^K -1\}\times \{0,1,\dots ,2^K
-1\}\\
&({}^{*}\!S \Rightarrow {}^{*}\!T)]\quad \wedge \\
&[\exists N \in {}^{*}\! \mathbb{N}\forall J>N\forall K\geq J
\exists X\subseteq \{0,1,\dots ,2^K -1\}\times \{0,1,\dots ,2^K
-1\}\\ & ({}^{*}\!S \Rightarrow {}^{*}\!T)],
\end{array} \] where ${}^{*}\!S$ and ${}^{*}\!T$ are the transferred versions of
the statements $S$ and $T$. Note that this necessitates replacing
only $A$ with ${}\nsa$ in the original $S$ and $T$.

We now choose any sufficiently large $J\in \nsn \setminus
\mathbb{N}$. The statement will still hold if we set $K=J$. This
results in a ``cover" of ${}\nsa$ by intervals of diameter $2^{-K}$.
Set \[ \at = \left\{\frac{j}{2^K} : (j-1,j) \in X\right\},\] where
$X$ is the set the existence of which is guaranteed in the second
line of the previous transferred statement.

By the transferred statement we know that $\sum_{(i,j)\in X} \left(
\frac{j-i}{2^K}\right)^{\beta} >1$, but $j-i =1$ because of the
choice of $K$ --- all the infinitesimal intervals are now of the
same size. Also, $\AT = |X|$; therefore $\frac{\AT}{2^{K \beta}}
>1$. Thus,
\[\textrm{meas}_{\beta}A>0 \Rightarrow \exists A_{\mathbf{T}}\subseteq \mathbf{T}, K\in {}\nsn
\setminus \mathbb{N} \textrm{ such that } \stan A_{\mathbf{T}}= A
\textrm{ and } \stan \left(\frac{\AT}{2^{K\beta}}\right) >0.\]

Since the converse holds by the nonstandard Frostman lemma, the
theorem is proved.\hfill{$\square$}

We now show that any set which satisfies certain of the above
properties is rich enough to yield Hausdorff dimension.

\begin{thm}
Consider a hyperfinite time line $\mathbf{T}$ based on the
infinitesimal $2^{-N}$, for a given $N\in \nsn\setminus \mathbb{N}$.
Suppose that an internal subset $A'$ of the time line is such that
$\stan( A') = A$ and for some $\alpha >0$
\begin{eqnarray}
\stan \left( \frac{|A'|}{2^{N\beta}}\right) &>& 0 \textrm{ for }
\beta
< \alpha \textrm{ and }\\
\stan \left( \frac{|A'|}{2^{N\beta}}\right) &=& 0 \textrm{ for }
\beta
> \alpha.
\end{eqnarray}
Then $\alpha = \textrm{dim} A$.
\end{thm}
{\bf{Proof.}} Given (4.1), the nonstandard version of Frostman's
lemma immediately implies that $\textrm{dim} A \geq \alpha$. For the
converse inequality, notice that the second condition implies that
for each $\varepsilon \in \mathbb{R}$, $\varepsilon >0$,
\[\frac{|A'|}{2^{N\beta}}< \varepsilon,\] which implies the
following nonstandard statement for each positive $\varepsilon \in
\mathbb{R}$:
\[ \begin{array}{lr} &\exists N \in {}\nsn \exists Y\subseteq \{0,1,\dots , 2^N -1\}
\forall x\in A' \exists i\ \in Y \left( x\in
(i2^{-N},(i+1)2^{-N}]\right) \wedge \\  & \left(
\frac{|Y|}{2^{N\beta}}< \varepsilon \right).\end{array}\]

The statement merely affirms the existence of an indexing set $Y$
for intervals of length $2^{-N}$ which form a cover of $A'$ and for
which the term $|Y|2^{-N\beta}$ is smaller than any real number.

Transferring down to the standard case, we find that for each
$\varepsilon >0$,
\[\begin{array}{lr} &\exists n\in \mathbb{N}\exists y\subseteq \{0,1,\dots , 2^n -1\}
 \forall x\in {} \stan\! A' \exists i\in y \left( x\in
(i2^{-n},(i+1)2^{-n}]\right) \wedge \\ & \left(
\frac{|y|}{2^{n\beta}} < \varepsilon \right) .\end{array}\] This
implies that $\textrm{meas}_{\beta}A =0$ and therefore that
$\textrm{dim}A \leq \alpha$.$\hfill{\square}$

For computational purposes it is therefore enough to find a set in
the time line with standard part $A$ that satisfies the conditions
in the above theorem. This fact will be used in subsequent sections.

In the sequel I refer to $\AT  \dt^{\beta}$ (where $\dt = 1/N$) as
nonstandard $\beta$-Hausdorff measure and to $\measb{A}$ as just
$\beta$-Hausdorff measure.

Several of the properties of the standard $\beta$-Hausdorff measure
can easily be seen to be valid in the nonstandard case, such as its
outer measure properties, invariance under translation (and
rotation, in the multidimensional case) and homogeneity of degree
$\beta$ with respect to dilation.

To illustrate some applications of this formulation, I first turn to
the perennial example of a set of non-integer dimension, the triadic
Cantor set.
 The ``base-infinitesimal" of the construction is\\
$\filt{1 , 3^{-1}, 3^{-2}, \dots , 3^{-k} , \dots } = \dt=1/N$. The
cardinality of the NS Cantor set $\AT$ is given by\\
$\filt{1,2/3,4/9,\dots ,(2/3)^k ,\dots } N$. The NS
$\beta$-Hausdorff measure of $A$ is then given by
\begin{eqnarray*}
\AT  \dt^{\beta} &=& \filt{(2/3)^k} N \filt{(1/3)^{k \beta}}\\
&=&\filt{(2/3^{\beta})^k },
\end{eqnarray*}
where I have used the obvious notation, $\filt{a^k }$ instead of
\mbox{$\filt{a, a^2 ,\dots , a^k, \dots}$.} The above expression
then has value $1$ for $\beta = \log 2/\log 3$, which is then $\DM
A$ by our previous theorems. Since the standard $\beta$-Hausdorff
sum for the triadic Cantor set is also $1$ for $\beta = \DM A$, I
suspect that the standard parts of the nonstandard sum will be equal
to the standard sum at $\DM A$ for other sets as well. This remains
to be proved.

\setcounter{section}{3}
\section{THE FRACTAL GEOMETRY OF BROWNIAN MOTION}
In this section I briefly discuss a nonstandard version of Brownian
local time, level sets and the effect of a Brownian motion on a set
with a given dimension. Although these results are not original, the
proofs using a nonstandard version of Hausdorff dimension are very
simple and intuitive. We start with a discussion on Brownian motion
in the nonstandard context, with emphasis on Anderson's simple and
beautiful construction~\cite{Anderson}.
\subsection{Anderson's construction of Brownian motion}
The idea is to construct Brownian motion as a hyperfinite random
walk, instead of, as is often done, a limit of random walks. We
start with a hyperfinite time line $\mathbf{T}$, based on a fixed $N
\in {}\nsn \setminus \mathbb{N}$. We let $\Omega
=\{-1,+1\}^{\mathbf{T}}$. If $\omega \in \Omega$, we define the
hyperfinite random walk as a polygonal path, filled in linearly
between time points $t\in \mathbf{T}$ with $B(\omega ,0)=0$ and
\[B(\omega ,t+\dt)-B(\omega ,t)=\vartriangle \! B(t)= \omega(s) \sqrt{\dt},\]
where $\omega (s)=\pm 1$. We let $\mathcal{C}_N$ be the set of all
such paths, $\mathcal{A}_N = {}^*\! \mathcal{P(C})_N$ and $W_N$ the
counting probability on $\mathcal{C}_N$ (where $\mathcal{P(C)}_N$
denotes the power set of $\mathcal{C}$ . This gives us the internal
probability space
\[(\mathcal{C}_N, \mathcal{A}_N, W_N)\]
which in turn gives us the Loeb space
\[\mathbf{\Omega} = (\mathcal{C}_N, L(\mathcal{A}_N), P_N = (W_N)_L).\]
The following theorem is due to Anderson~\cite{Cutland}. Recall that
an internal function $F$ is S-continuous if, whenever arguments $x$
and $y$ are infinitesimally close, the corresponding function values
$F(x)$ and $F(y)$ are infinitesimally close as well.

\begin{thm}
$1.$ For almost all $B\in \mathcal{C}_N$, $B$ is S-continuous and
gives a continuous path $b={}\stan B \in \mathcal{C}$.

$2.$ For Borel $D\subseteq \mathcal{C}$, $ W(D) = P_N (st^{-1}(D)) $
is Wiener measure.

$3.$ The following process is a Brownian motion on the space
$\mathbf{\Omega}$:
\[b(t, \omega ) = {}\stan B(w,t):[0,1]\times \Omega \to \mathbb{R}.\]
\end{thm}
For a proof, as well as a nonstandard version of the central limit
theorem, see~\cite{Albeverio}.
\subsection{Brownian local time}
The local time of a Brownian motion is a measure of the time a
Brownian motion spends at $x$, giving an indication of how many
times the path returns to a certain value. The Lebesgue measure of
this set is $0$, but it can be described using Hausdorff measure, as
we shall see shortly.

Heuristically, we define the local time $l(t,x)$ as
\[l(t,x) = \int_{0}^{t} \delta (x-b(s))ds,\]
where $b$ is a Brownian motion and $\delta$ the delta function. A
more precise definition and a discussion of local time can be found
in~\cite{ItoMcKean}. The integral therefore ``counts" how many times
the Brownian path visits $x$ up until the time $t$. The standard
approach (which can be found in detail in, for example,
~\cite{ItoMcKean}) is to show there exists a jointly continuous
process $l(t,x)$ such that
\[l(t,x) = \frac{d}{dx} \int_{0}^{t} I_{(-\infty ,x]}(b(s)) ds,\]
for almost all $(t,x)\in [0,1]\times \mathbb{R}$, where $I_A$ is the
characteristic function of the set $A$. Note that although the
definition is valid for a time line $[0,\infty)$ as well as $[0,1]$,
I use a bounded interval throughout. The nonstandard approach, due
to Perkins~\cite{Perkins}, is clearer and more intuitive. We think
of the Brownian path $b$ as the standard part of a hyperfinite
random walk. The following exposition follows~\cite{Albeverio}. We
start by approximating $l(t,x)$ by
\[ (\vartriangle x)^{-1} \int_{0}^{t} I_{[x,x+\vartriangle x]}
(b(s))ds.\] Now replace the time line $[0,1]$ by a discrete
hyperfinite time line $\mathbf{T}$ and the space $\mathbb{R}$ by
\\ $\Gamma = \{0, \pm \sqrt{\dt}, \dots, \pm n\sqrt{\dt}, \dots, \pm
N\sqrt{\dt}\}$ and define the internal process $L:\mathbf{T} \times
\Gamma \to {}^*\mathbb{R}$ by
\[L(t,x)= \sum_{s<t} I_{x}(B(s))(\dt)^{1/2},\]
where $I_x = I_{\{x\}}$. Perkins showed that $L(t,x)$ has a standard
part which is Brownian local time. He used the nonstandard
formulation to prove the following global characterisation of local
time, which was previously known to hold only for each $x$
separately: Let $\lambda (t,x,\delta)$ be the Lebesgue measure of
the set of points within a distance of $\delta /2$ of $\{s\leq
t|b(s)=x\}$. Then for almost all $\omega \in \Omega$ and each $t_0 >
0$,
\[\lim_{\delta \to 0^+} \sup_{t\leq t_0 ,x\in \mathbb{R}}
|\lambda(t,x,\delta )\delta^{-1/2} -2(2/\pi )^{1/2} l(t,x)|=0.\] It
is shown in~\cite{ItoMcKean} that local time for $t=1$ is the same
as $\frac{1}{2}$-dimensional Hausdorff measure on a level set of
Brownian motion; that is, for fixed $x$, $l(1,x) =
\textrm{meas}_{\frac{1}{2}}(b^{-1}(x))$. From the nonstandard
formulation, however, it is immediately clear. If we define the set
$A$ as the set of all $t\in [0,1]$ such that $b(t)=x$, the
nonstandard local time becomes simply $|A_{\mathbf{T}}|\dt^{1/2}$,
where $A_{\mathbf{T}}$ is the nonstandard version of the set $A$
encountered in the proof of Theorem 3.3. But this is exactly the
quantity whose standard part is the same as
$\frac{1}{2}$-dimensional Hausdorff measure (up to a finite constant
factor --- which depends on which author one reads). We must now
show that level sets have dimension $1/2$. We show this for $x=0$
only, since the level sets all have similar dimension. Denote the
zero set of a Brownian path $b(\omega)$ by $Z_{\omega}$ (or just $Z$
when possible). The required nonstandard version of the set is
denoted by $Z_{\omega ,\mathbf{T}}$. We now turn to a standard
property of local time to show that the dimension of this set is
$1/2$. It can be shown (as for instance in~\cite{ItoMcKean}) that
local time is identical in law to the function
\[M_{\omega}(t) = \max_{s\leq t}b(s).\]
This implies that $P[l(1,0) >0]=P[\max_{s\leq 1}b(s)>0]=1$. By the
nonstandard formulation of local time, this immediately implies that
$\stan (Z_{\omega, \mathbf{T}}) >0$, which implies that $\DM Z \leq
1/2$, almost surely. By the same token, however, $l(1,0)$ is almost
certainly finite, implying that $\stan (Z_{\omega, \mathbf{T}}) <
\infty$ and therefore $\DM(Z) \geq 1/2$.

The following lemma will be used in the subsequent section. In this
case the standard approach is easier than the hyperfinite, by using
the H\"{o}lder condition for Brownian motion. The proof is akin to
the proof of the dimension of the level set found in~\cite{Peres}.
We will use the fact that $Y(t)=M(t)-b(t)$ (where $M$ is as defined
above and $b$ is a standard Brownian motion) has the same
distribution as $b$~\cite{Hida}. Note also that the zeroes of $Y$
correspond to the global maxima (from the left) of $b$; these are
known as \emph{record times}.

\begin{lem}
If $Z$ is a level set and $A\subseteq Z$, then either $D$ has
dimension $1/2$ or $\textrm{dim}Z/A = 1/2$.
\end{lem}

{\bf{Proof.}} Since $M(t)$ is an increasing function, it can be
considered to be the distribution of some measure $\mu$, which has
its support on the set of record times. Let $Z$ be the zeroset for
$Y$. (Because of the similar distributions, dimensional results for
this set will hold for any Brownian level set.) We therefore have
that $\mu(a,b] = M(b)-M(a)$. By the H\"{o}lder condition for
Brownian motion we get
\[M(b) -M(a) \leq \max_{0\leq h \leq b-a} b(a+h)-b(a) \leq
C_{\alpha}(b-a)^{\alpha}\] for some constant $C_{\alpha}$, for all
$\alpha<1/2$.

Consider now a subset $A$ of $Z$ and let $\mu'$ be the restriction
of the measure $\mu$ to $A$; that is,
\[\mu '(a,b] = \mu \{(a,b]\cap A\}.\]
Suppose that $\mu'$ is not $0$ on every interval. Then,
\[0< \mu '(a,b]\leq M(b)-M(a)\leq C_{\alpha}(b-a),\]
as above, for some $(a,b]$. By normalising the measure $\mu '$ we
find a measure $\nu$ such that $\nu (a,b] \leq D_{\alpha}(b-a)$ for
some $D_{\alpha}$, for each $0<\alpha <1/2$. By Frostman's lemma,
$\textrm{dim}A\geq 1/2$.

If there exists no interval $(a,b]$ such that $\mu '(a,b] >0$, then
we can safely leave out $A$ without changing the Hausdorff dimension
of $Z$, as in~\cite{Peres}, and $\textrm{dim}Z/A
=1/2$.\hfill{$\square$}

\begin{cor} If $\dim A<1/2$, then the inverse image of any element
in $B(A)$ (where $B$ is a Brownian motion) has dimension $0$.
\end{cor}
\subsection{The image of a set under Brownian motion}
A very interesting property of Brownian motion is its effect on sets
of a certain Hausdorff dimension. If a compact subset of $[0,1]$ has
dimension $\alpha < 1/2$, its image under Brownian motion is a set
of dimension $2\alpha$. A set of dimension $\alpha >1/2$ will have
dimension $1$ and will almost surely contain an interval. As for
sets of dimension $1/2$, we have seen above that they may have an
image of dimension $0$. No hard and fast rule exists for such sets.
We now look at nonstandard proofs of these results. The advantage of
this approach is a more intuitive (counting) argument. A Fourier
analytical proof of the following can be found in~\cite{Kahane}.
\begin{thm} Let $A\subset [0,1]$ be a compact set. If $\DM A
= \alpha <1/2$ and $b$ is a Brownian motion, $dim\, b(A)=2\alpha$.
\end{thm}
\textbf{Proof.} The basis for the time line of the image is no
longer $\dt$, but $\sqrt{\dt}$. Let $B$ denote a nonstandard
Brownian motion which has $b$ as standard part. We let
$A_{\mathbf{T}}$ be the nonstandard counterpart of $A$ constructed
in Theorem 3.3. Since $|B(A_T )|\leq |A_T |$ and we know that $|A_T
|\dt^{\beta} \approx 0$ for $\beta > \alpha$, we will have that
$|B(A_T )|\dt^{\beta}\approx 0$ for any $\beta
>\alpha $. Therefore, $|B(A_T )|(\sqrt{\dt})^{\gamma}\approx 0$
for $\gamma >2\alpha$ and we conclude that $\textrm{dim}b(A)\leq
2\textrm{dim}A$ by Theorem 3.4, since $ ^{\circ}\!B(A_{\mathbf{T}})=
b(A)$ (because of the S-continuity of the functions involved). It is
left to show that $\textrm{dim}b(A)\geq 2\textrm{dim}A$. This is not
quite as simple as the previous proof, since the matter of possible
level sets complicates the question of the cardinality of the image.
We overcome this by considering only one element of each level set
and discarding the rest. The remaining set will have the same
dimension as the original and the image will have the same
cardinality. This is made possible because the set $A$ has a
dimension of less than $1/2$. Any subsets of level sets in $A$ are
small enough to be left out (mostly) without affecting the dimension
(see Lemma 4.2). For any $x\in b(A)$, pick one representative $x'
\in B(A_{\mathbf{T}})$ such that $\stan x' = x$ and some $t \in
A_{\mathbf{T}}$ such that $B(t)=x'$. Let $X'$ be the set of all such
representatives $x'$. We denote by $L_{x',\mathbf{T}}$ the subset of
$A_{\mathbf{T}}$ for which $\stan B(L_{x',\mathbf{T}}) =x$.We want
to show now that the standard parts of the sums
\[\sum_{x' \in X'} \frac{1}{N^{\alpha}},\quad
\sum_{x' \in X'} \frac{|L_{x',\mathbf{T}}|}{N^{\alpha}}\] are $0$
and $\infty$ for the same values of $\alpha$. To do this, all that
is necessary is to show that the first one is infinite whenever the
second one is. So suppose that
\[\stan \left( \sum_{x' \in X'}
\frac{|L_{x',\mathbf{T}}|}{N^{\alpha}} \right)= \infty.\] We know
that \[\frac{|L_{x',\mathbf{T}}|}{N^{\alpha}} = s_{x}^{\beta}
\approx 0\] for any $\beta >0$. This implies that
\[\sum_{x' \in X'} \frac{s_{x}^{\beta}
N^{\beta}}{N^{\alpha}}= \sum_{x' \in X'}
\frac{s_{x}^{\beta}}{N^{\alpha -\beta}}\leq \sum_{x'\in X'}
\frac{1}{N^{\alpha -\beta}}=\infty,\] for $\beta$ arbitrarily close
to $0$. Thus we may conclude that the number of level sets is
important and not the cardinality of each. But the number of level
sets is equal to the cardinality of the range, thus the standard
parts of \[ \frac{|B(A_T )|}{N^{\alpha}} \textrm{ and } \frac{|A_T
|}{N^{\alpha}}\] are $0$ and $\infty$ for the same values of
$\alpha$. Keeping in mind that the time line of the image is based
on $\sqrt{\dt}$ and not $\dt$, we can conclude that the dimension
has doubled.\hfill{$\square$}

%The nonstandard formulation of Hausdorff dimension seems to yield a
%more intuitive proof in some cases. This could have implications for
%the study of phenomena which have a satisfactory nonstandard
%version, such as the Anderson construction of Brownian motion. The
%work done so far has lead to dimensional results in the case of
%complex oscillations (a constructive version of Brownian motion -
%see~\cite{Fouche}, for instance), to be published in an upcoming
%paper.

As has been shown, this formulation of Hausdorff dimension proves to
be useful in studying the fractal properties of Brownian motion,
according to Anderson's construction. One may speculate that it
might also aid in other contexts where the phenomenon in question
has a satisfactory formulation in terms of the hyperfinite time
line. Specifically, the approach given here is based on a {\it
counting measure} approach and may therefore apply to systems where
the measure used may be constructed as a Loeb measure. (For more
interesting applications of Loeb measure to stochastic fluid
mechanics, stochastic calculus of variations and mathematical
finance theory, see~\cite{Cutland}.) The method has already been
used in the study of the points of rapid growth of Brownian motion,
as well as those of so-called complex oscillations~\cite{Fouche}, a
constructive version of Brownian motion, to be presented in a
forthcoming paper by the author.

%%%%%%%%%%%%%%%%%%%%%%%%%%%%%%%%%%%%%%%%%%%%%%%%%%%%%%%%%%%%%%%%%%%%%%%%%%%%%%%%%

\end{document}